\documentclass[reqno, 11pt]{amsart}
\usepackage[T1]{fontenc}
\usepackage{lmodern}
\usepackage[utf8]{inputenc}
\usepackage{geometry}
\usepackage{amssymb}
\usepackage{amsmath}
\usepackage{amsfonts}
\usepackage{indentfirst}
\usepackage{mathrsfs}
\usepackage[colorlinks=true,citecolor=blue,linkcolor=blue,urlcolor=blue]{hyperref}

\newtheorem{theorem}{Theorem}[section]
\newtheorem{lemma}[theorem]{Lemma}

\theoremstyle{definition}
\newtheorem{definition}[theorem]{Definition}
\theoremstyle{remark}

\numberwithin{equation}{section}
\newcommand{\R}{\mathbb R}

\begin{document}
\title[Upper-critical fractional Choquard equations]
{Normalized Solutions and Semiclassical Concentration for Upper-Critical Fractional Choquard Equations}

\author{Yergen Aikyn}
\address{Yergen Aikyn:
\endgraf
Department of Mathematics: Analysis, Logic and Discrete Mathematics, Ghent University, Belgium}
\email{aikynyergen@gmail.com}

\author{Yongpeng Chen}
\address{Yongpeng Chen:
\endgraf
School of Science, Guangxi University of Science and Technology, Liuzhou, China}
\email{yongpengchen@mail.bnu.edu.cn}

\author{Michael Ruzhansky}
\address{Michael Ruzhansky:
\endgraf
Department of Mathematics: Analysis, Logic and Discrete Mathematics, Ghent University, Belgium
\endgraf
School of Mathematical Sciences, Queen Mary University of London, United Kingdom}
\email{michael.ruzhansky@ugent.be}

\author{Zhipeng Yang}
\address{Zhipeng Yang:
\endgraf
Department of Mathematics, Yunnan Normal University, Kunming, China
\endgraf
Department of Mathematics: Analysis, Logic and Discrete Mathematics, Ghent University, Belgium}
\email{yangzhipeng326@163.com}

\thanks{Corresponding author: yangzhipeng326@163.com.}

\subjclass[2020]{35R11, 35J60, 35B40}
\keywords{Fractional Choquard equations; normalized solutions; semiclassical concentration}
\date{}

\begin{abstract}
We study a fractional Choquard equation with an upper-critical Hartree term,
an $L^2$-supercritical Hartree perturbation, and a semiclassical potential
under a prescribed $L^2$-mass constraint. The potential is bounded and
nonnegative, has a nonempty zero set, and has a positive lower limit at
infinity. For every prescribed mass and all sufficiently small semiclassical
parameters, we prove the existence of a pair of normalized solutions
$\pm u_\varepsilon$ with a negative Lagrange multiplier. The proof combines a
strict energy bound below the critical one-bubble level, compactness modulo
translations for the autonomous ground-state set, a simultaneous cutoff of
both Hartree terms, and a localized constrained mountain-pass argument.
Moreover, suitable translates of $u_\varepsilon$ converge strongly in
$H^s(\R^N)$ to a positive autonomous ground state, and the corresponding
concentration points approach the zero set of the potential as
$\varepsilon\to0$.
\end{abstract}

\maketitle

\section{Introduction and main results}

We consider the normalized fractional Choquard equation
\begin{equation}\label{eq1.1}
\begin{split}
(-\Delta)^s u+V(\varepsilon x)u
&=\lambda u+\bigl(I_\alpha*|u|^p\bigr)|u|^{p-2}u+\bigl(I_\alpha*|u|^q\bigr)|u|^{q-2}u
\quad\text{in }\R^N,
\end{split}
\end{equation}
subject to
\[
 \int_{\R^N}|u|^2\,dx=a.
\]
Here $a>0$ is prescribed, $\lambda\in\R$ is the associated Lagrange
multiplier, and $\varepsilon>0$ is the semiclassical parameter. We assume
\[
 N\ge2,\qquad 0<s<\min\left\{1,\frac N4\right\},\qquad
 N-4s<\alpha<N,
\]
and
\[
 \frac{N+2s+\alpha}{N}<q<p=2^*_{\alpha,s}:=
 \frac{N+\alpha}{N-2s}.
\]
The Riesz potential is
\[
 I_\alpha(x)=A_{N,\alpha}|x|^{\alpha-N},\qquad
 A_{N,\alpha}=
 \frac{\Gamma((N-\alpha)/2)}{\pi^{N/2}2^\alpha\Gamma(\alpha/2)}.
\]
The potential satisfies
\begin{equation}\tag{V}\label{eq1.2}
\begin{gathered}
 V\in C(\R^N)\cap L^\infty(\R^N),\qquad V\ge0,\qquad V(0)=0,\\
 V_\infty:=\liminf_{|x|\to\infty}V(x)\in(0,\infty).
\end{gathered}
\end{equation}
We write
\[
 M:=\{x\in\R^N:V(x)=0\}.
\]

For
\[
 B_r(u):=\int_{\R^N}(I_\alpha*|u|^r)|u|^r\,dx,
\]
the Hardy--Littlewood--Sobolev inequality \cite{MR1817225} and the
fractional Sobolev embedding imply that $B_r$ is well defined for
\[
 \frac{N+\alpha}{N}\le r\le2^*_{\alpha,s}.
\]
The mass-preserving dilation $u_t(x)=t^{N/2}u(tx)$ satisfies
\[
 \|(-\Delta)^{s/2}u_t\|_2^2
 =t^{2s}\|(-\Delta)^{s/2}u\|_2^2,
 \qquad
 B_r(u_t)=t^{Nr-(N+\alpha)}B_r(u).
\]
Hence the $L^2$-critical Hartree exponent is
\[
 \bar p=\frac{N+2s+\alpha}{N}.
\]
Both nonlinearities in \eqref{eq1.1} are therefore $L^2$-supercritical,
whereas the $p$-term has HLS upper-critical growth. The constrained energy
is unbounded from below, and the critical Hartree term introduces the usual
loss of compactness.

Choquard equations arise in the one-component plasma model analyzed by Lieb
\cite{Lieb1977}, in Pekar's polaron model \cite{Pekar1954}, and in
Schr\"odinger--Newton models \cite{Diosi1984,MorozPenroseTod1998}. General
existence, regularity, symmetry, and decay results for the classical equation
are developed in \cite{MorozVanSchaftingen2013,MR3625092}. Fractional models
are related to Laskin's fractional quantum mechanics \cite{Laskin2000}.
Basic properties of the fractional Laplacian and fractional Sobolev spaces
are reviewed in \cite{MR2944369}.

Normalized solutions are constrained critical points on an $L^2$-sphere.
The variational framework was introduced for nonlinear Schr\"odinger equations
by Jeanjean \cite{Jeanjean1997NA} and subsequently developed in
\cite{BartschDeValeriola2013ARMA,Soave2020JDE,Soave2020JFA}. For fractional
Choquard equations, Li and Luo \cite{LiLuo2020} treated the
$L^2$-supercritical, HLS-subcritical regime. The HLS lower-critical case and
its nonautonomous counterpart were studied in
\cite{YuTangZhang2023,LiuSunZhang2024}. At the HLS upper-critical exponent,
normalized solutions were obtained in
\cite{LanHeMeng2023,HeRadulescuZou2022,MengHe2024} by combining sharp
inequalities, truncation, and concentration--compactness.

Mixed Hartree nonlinearities and external potentials have been considered in
several recent works. A critical fractional $p$-Laplacian problem was studied
in \cite{ZhangThinLiang2025}. Weighted fractional Choquard equations with
mixed nonlinearities up to the $L^2$-critical range were treated in
\cite{ChenYangZhangPotentials2025}, while the corresponding unweighted mixed
problem was studied in \cite{MR4921761}. In the local Choquard setting,
Meng, He, and Winkert \cite{MengHeWinkert2026} proved multiplicity and
concentration for HLS-subcritical growth.

The present paper addresses the fully $L^2$-supercritical range
\[
 \bar p<q<p=2^*_{\alpha,s}
\]
for arbitrary prescribed mass. The autonomous analysis gives a positive
ground level strictly below the energy of one critical Hartree bubble,
strict monotonicity with respect to the mass, and compactness of ground states
modulo translations. For the nonautonomous problem, a single cutoff is
applied simultaneously to the two Hartree terms. The constrained deformation
is localized near the positive autonomous ground-state orbit, where the
cutoff is inactive and the Lagrange multipliers remain uniformly negative.
This yields a normalized critical point and identifies its concentration
profile as $\varepsilon\to0$.

Our main result is the following.

\begin{theorem}\label{Thm1.1}
Assume that $N\ge2$, $0<s<\min\{1,N/4\}$, $N-4s<\alpha<N$,
\[
 \frac{N+2s+\alpha}{N}<q<p=\frac{N+\alpha}{N-2s},
\]
and that \textnormal{(V)} holds. For every $a>0$, there exists
$\varepsilon_a>0$ such that, for each $0<\varepsilon<\varepsilon_a$, there
are $u_\varepsilon\in H^s(\R^N)$ and $\lambda_\varepsilon<0$ for which
\[
 (u_\varepsilon,\lambda_\varepsilon),\qquad
 (-u_\varepsilon,\lambda_\varepsilon)
\]
are distinct weak solutions of \eqref{eq1.1} and
\[
 \|u_\varepsilon\|_2^2=a.
\]
Moreover, if $\varepsilon_n\to0$, then, after passing to a subsequence,
there exist $z_n\in\R^N$ and a positive ground state $w$ of the autonomous
zero-potential problem such that
\[
 u_{\varepsilon_n}(\,\cdot+z_n)\longrightarrow w
 \quad\text{in }H^s(\R^N),
 \qquad
 \operatorname{dist}(\varepsilon_n z_n,M)\longrightarrow0.
\]
\end{theorem}

Section~\ref{Sec2} contains the preliminary results. Section~\ref{Sec3}
studies the autonomous ground level and compactness modulo translations.
Section~\ref{Sec4} establishes the cutoff and the local Palais--Smale theory,
and Section~\ref{Sec5} proves Theorem~\ref{Thm1.1}.

\section{Preliminaries}\label{Sec2}

We work with real-valued functions. The symbols \(C\) and \(C_i\) denote
positive constants that may vary from line to line, and \(o_n(1)\to0\) as
\(n\to\infty\). We fix the functional setting and recall the
auxiliary inequalities used below. For \(s\in(0,1)\), the fractional
Sobolev space \(H^s(\R^N)\) has the equivalent characterizations
\[
\begin{aligned}
H^s(\R^N)
&=\Big\{u\in L^2(\R^N):\ \frac{u(x)-u(y)}{|x-y|^{\frac{N}{2}+s}}\in L^2(\R^N\times\R^N)\Big\} \\
&=\Big\{u\in L^2(\R^N):\ \int_{\R^N}\big(1+|\xi|^{2s}\big)|\mathcal{F}(u)(\xi)|^2\,d\xi<\infty\Big\},
\end{aligned}
\]
where $\mathcal{F}(u)$ denotes the Fourier transform of $u$.

The Gagliardo seminorm of $u$ is given by
\[
[u]_{H^s(\R^N)}^2
=\int_{\R^N}\int_{\R^N}
\frac{|u(x)-u(y)|^2}{|x-y|^{N+2s}}\,dx\,dy,
\]
and a standard computation using \cite[Propositions 3.4 and 3.6]{MR2944369} yields, for $u\in H^s(\R^N)$,
\[
\int_{\R^N}\big|(-\Delta)^{s/2}u\big|^2\,dx
=\int_{\R^N}|\xi|^{2s}|\mathcal{F}(u)(\xi)|^2\,d\xi
=\frac{1}{2}C(N,s)\int_{\R^N}\int_{\R^N}
\frac{|u(x)-u(y)|^2}{|x-y|^{N+2s}}\,dx\,dy,
\]
where $C(N,s)>0$ is the normalizing constant in the definition of $(-\Delta)^s$.

Throughout the paper we endow $H^s(\R^N)$ with the norm
\[
\|u\|_{H^s(\R^N)}
=\Big(\int_{\R^N}|u|^2\,dx
+\int_{\R^N}\big|(-\Delta)^{s/2}u\big|^2\,dx\Big)^{\frac{1}{2}},
\]
and with the associated inner product
\[
\langle u,v\rangle
=\int_{\R^N}(-\Delta)^{s/2}u\,(-\Delta)^{s/2}v\,dx
+\int_{\R^N}u v\,dx.
\]
We use the homogeneous fractional Sobolev space
\[
\mathcal{D}^{s,2}(\R^N)
=\Big\{u\in L^{2_s^*}(\R^N):
\int_{\R^N}\int_{\R^N}
\frac{|u(x)-u(y)|^2}{|x-y|^{N+2s}}\,dx\,dy<\infty\Big\},
\]
equipped with the norm
\[
\|u\|
=\Big(\int_{\R^N}\big|(-\Delta)^{s/2}u\big|^2\,dx\Big)^{\frac{1}{2}},
\]
where $2_s^*=\frac{2N}{N-2s}$ is the fractional critical Sobolev exponent.

We recall the sharp Hartree interpolation inequality from
\cite{feng-zhang}.

\begin{lemma}\label{Lem2.1}
Assume the standing conditions and let $t$ satisfy
\[
\frac{N+\alpha}{N}<t<2^*_{\alpha,s},
\qquad
2^*_{\alpha,s}=\frac{N+\alpha}{N-2s}.
\]
Then there exists $C_{\alpha,t}>0$ such that, for all $u\in H^s(\R^N)$,
\[
\int_{\R^N}\big(I_\alpha*|u|^t\big)|u|^t\,dx
\le C_{\alpha,t}\,\|u\|^{2 t \gamma_{t,s}}\|u\|_2^{2 t(1-\gamma_{t,s})},
\]
where
\[
\gamma_{t,s}=\frac{N t-N-\alpha}{2 s t}\in(0,1),
\]
and the optimal constant $C_{\alpha,t}$ is given by
\[
C_{\alpha,t}
=\frac{2 s t}{2 s t-N t+N+\alpha}
\Big(\frac{2 s t-N t+N+\alpha}{N t-N-\alpha}\Big)^{\frac{Nt-N-\alpha}{2 s}}\,
\|U\|_2^{2-2 t},
\]
where $U$ is a ground-state solution of
\[
(-\Delta)^s U+U-\big(I_\alpha *|U|^t\big)|U|^{t-2} U=0
\quad\text{in }\R^N.
\]
\end{lemma}

At the HLS upper critical exponent $2^*_{\alpha,s}$, the optimal constant
is determined by the critical Choquard inequality
\cite{he-r}.

\begin{lemma}\label{Lem2.2}
Assume the standing conditions and define
\[
\mathcal{S}_\alpha
=\inf_{u\in\mathcal{D}^{s,2}(\R^N)\setminus\{0\}}
\frac{\displaystyle\int_{\R^N}\big|(-\Delta)^{s/2}u\big|^2\,dx}
{\displaystyle\Big(\int_{\R^N}\big(I_\alpha*|u|^{2^*_{\alpha,s}}\big)|u|^{2^*_{\alpha,s}}\,dx\Big)^{\frac{1}{2^*_{\alpha,s}}}}.
\]
Then a function
$u\in\mathcal{D}^{s,2}(\R^N)\setminus\{0\}$
attains $\mathcal{S}_\alpha$ if and only if
\[
u(x)=c\Big(\frac{\varepsilon}{\varepsilon^2+|x-x_0|^2}\Big)^{\frac{N-2s}{2}},
\qquad x\in\R^N,
\]
for some $x_0\in\R^N$, $c\in\R\setminus\{0\}$, and
$\varepsilon>0$.
\end{lemma}

With $\gamma_{2^*_{\alpha,s}, s}:=1$ and
\[
C_{\alpha,2^*_{\alpha,s}}
=\big(\mathcal{S}_\alpha^{-1}\big)^{2^*_{\alpha,s}},
\]
Lemma~\ref{Lem2.1} extends to every
\[
t\in\Big(\frac{N+\alpha}{N},\,2^*_{\alpha,s}\Big].
\]

The following lemma gives the local compactness of the Riesz potential used
below.

\begin{lemma}\label{Lem2.3}
Set
\[
 t=\frac{2N}{N+\alpha},\qquad t'=\frac{2N}{N-\alpha}.
\]
Let $(F_n)$ be bounded in $L^t(\R^N)$ and satisfy $F_n\to0$ almost
everywhere. Then, for every bounded measurable set $K\subset\R^N$
and every $1\le k<t'$,
\[
 I_\alpha*F_n\longrightarrow0\quad\text{in }L^k(K).
\]
Moreover, let $d>1$ and assume that
\[
 \frac1k=\frac1{t'}+\frac1d<1.
\]
If $G\in L^d(\R^N)$, then
\[
 (I_\alpha*F_n)G\longrightarrow0\quad\text{in }L^k(\R^N).
\]
If $F\in L^t(\R^N)$ and $(G_n)$ is bounded in $L^d(\R^N)$
with $G_n\to0$ almost everywhere, then
\[
 (I_\alpha*F)G_n\longrightarrow0\quad\text{in }L^k(\R^N).
\]
\end{lemma}

\begin{proof}
Boundedness and almost-everywhere convergence imply
$F_n\rightharpoonup0$ in $L^t$. Fix $k<t'$ and choose
$\widehat k\in(\max\{k,t\},t')$. Let $a>1$ be determined by
\[
 1+\frac1{\widehat k}=\frac1t+\frac1a.
\]
Then $a<N/(N-\alpha)$. Decompose the Riesz kernel as
\[
 I_\alpha=I_\alpha\mathbf1_{B_\delta}
 +I_\alpha\mathbf1_{B_R\setminus B_\delta}
 +I_\alpha\mathbf1_{\R^N\setminus B_R}.
\]
Young's inequality gives
\[
 \|(I_\alpha\mathbf1_{B_\delta})*F_n\|_{\widehat k}
 \le \|I_\alpha\mathbf1_{B_\delta}\|_a\|F_n\|_t,
\]
and the right-hand side tends to zero uniformly in $n$ as $\delta\downarrow0$.
Since $I_\alpha\mathbf1_{B_R\setminus B_\delta}\in L^{t'}$, weak
convergence gives pointwise convergence of the middle convolution, while
H\"older's inequality gives a uniform bound. Hence this term converges to
zero in $L^{\widehat k}(K)$. Finally,
\[
 \|(I_\alpha\mathbf1_{\R^N\setminus B_R})*F_n\|_\infty
 \le \|I_\alpha\mathbf1_{\R^N\setminus B_R}\|_{t'}\|F_n\|_t,
\]
and the right-hand side tends to zero as $R\to\infty$. Thus
$I_\alpha*F_n\to0$ in $L^{\widehat k}(K)$, and therefore in $L^k(K)$.

For the first product convergence, approximate $G$ in $L^d$ by bounded
compactly supported functions and use the local convergence just proved,
together with the HLS estimate in $L^{t'}$. For the second, approximate
$I_\alpha*F$ in $L^{t'}$ by bounded compactly supported functions. On a
fixed bounded set, $G_n\to0$ in $L^k$ because $k<d$, while the approximation
error is controlled by H\"older's inequality. This proves the lemma.
\end{proof}

\begin{lemma}\label{Lem2.4}
Let $W\in C(\R^N)\cap L^\infty(\R^N)$,
$y_n\to y$, and $\varepsilon_n\to0$. Then, for every
$\varphi\in L^2(\R^N)$,
\[
 \|[W(y_n+\varepsilon_n\,\cdot)-W(y)]\varphi\|_2\longrightarrow0.
\]
If $\mathcal K\subset L^2(\R^N)$ is compact, then
\[
 \sup_{\varphi\in\mathcal K}
 \int_{\R^N}|W(y_n+\varepsilon_n x)-W(y)|\,|\varphi(x)|^2\,dx
 \longrightarrow0.
\]
Furthermore, if $\varphi_n\to\varphi$ in $L^2(\R^N)$, then
\[
 \int_{\R^N}W(y_n+\varepsilon_n x)|\varphi_n(x)|^2\,dx
 \longrightarrow W(y)\|\varphi\|_2^2.
\]
\end{lemma}

\begin{proof}
Fix $\varphi\in L^2(\R^N)$ and $\delta>0$. Choose $R>0$ such that
\[
 \int_{\R^N\setminus B_R}|\varphi|^2\,dx<\delta.
\]
For all sufficiently large $n$,
\[
 \sup_{x\in B_R}|W(y_n+\varepsilon_n x)-W(y)|<\delta.
\]
Consequently,
\[
 \|[W(y_n+\varepsilon_n\,\cdot)-W(y)]\varphi\|_2^2
 \le \delta^2\|\varphi\|_2^2+4\|W\|_\infty^2\delta,
\]
which proves the first assertion.

Let $\mathcal K\subset L^2(\R^N)$ be compact. A finite $L^2$-net for
$\mathcal K$ and the preceding tail estimate give
\[
 \lim_{R\to\infty}\sup_{\varphi\in\mathcal K}
 \int_{\R^N\setminus B_R}|\varphi|^2\,dx=0.
\]
The preceding estimate is therefore uniform on $\mathcal K$, which
proves the second assertion. Finally,
\[
 \begin{aligned}
 &\left|\int_{\R^N}W(y_n+\varepsilon_n x)|\varphi_n|^2\,dx
      -W(y)\|\varphi\|_2^2\right|\\
 &\quad\le \|W\|_\infty
      \big\||\varphi_n|^2-|\varphi|^2\big\|_1
 +\int_{\R^N}|W(y_n+\varepsilon_n x)-W(y)|\,|\varphi|^2\,dx.
 \end{aligned}
\]
The first term tends to zero by $L^2$-convergence, and the second tends to
zero by the first assertion.
\end{proof}

\begin{lemma}\label{Lem2.5}
Let $m\in\mathbb N$, $\omega>0$, and
\[
 H_j\in L^{2N/\alpha}(\R^N)
       +L^{2N/(\alpha+2s)}(\R^N),
 \qquad 1\le j\le m.
\]
If $u\in H^s(\R^N)$ is a weak solution of
\[
 (-\Delta)^s u+\omega u
 =\sum_{j=1}^m(I_\alpha*(H_j u))H_j
 \quad\text{in }\R^N,
\]
then
\[
 u\in L^\ell(\R^N)
 \qquad\text{for every }2\le\ell<\frac N\alpha\,2_s^*.
\]
\end{lemma}

\begin{proof}
Set
\[
 H:=\sum_{j=1}^m|H_j|.
\]
Since the sum is finite,
\[
 H\in L^{2N/\alpha}(\R^N)+L^{2N/(\alpha+2s)}(\R^N).
\]
Positivity of the Riesz kernel gives, for $u,\varphi\in H^s(\R^N)$,
\[
 \left|\sum_{j=1}^m\int_{\R^N}
 (I_\alpha*(H_j u))H_j\varphi\,dx\right|
 \le\int_{\R^N}(I_\alpha*(H|u|))H|\varphi|\,dx.
\]
Let $2\le\ell<2N/\alpha$ and let $u_L$ be the truncation of $u$ at
height $L$. In the proof of
\cite[Proposition~3, equations (3.16)--(3.24)]{Ambrosio2025Remarks}, test
the equation by $|u_L|^{\ell-2}u_L$ and split the domain into
$\{|u|\le L\}$ and $\{|u|>L\}$. The contribution from the first set is
bounded by
\[
 \begin{aligned}
 &\sum_{j=1}^m\int_{\R^N}
 (I_\alpha*(|H_j||u_L|))|H_j||u_L|^{\ell-1}\,dx\\
 &\qquad\le
 \int_{\R^N}(I_\alpha*(H|u_L|))H|u_L|^{\ell-1}\,dx.
 \end{aligned}
\]
Apply \cite[Lemma~7]{Ambrosio2025Remarks} with
\[
 v=|u_L|^{\ell/2},\qquad \theta=\frac2\ell.
\]
The condition $\ell<2N/\alpha$ gives
$\alpha/N<\theta\le1<2-\alpha/N$, so the gradient term can be absorbed.
The contribution from $\{|u|>L\}$ tends to zero as $L\to\infty$ by the
Hardy--Littlewood--Sobolev inequality and dominated convergence, exactly as
in equations (3.23)--(3.24) of the cited proof. Letting $L\to\infty$ gives
\[
 \|u\|_{\ell 2_s^*/2}^{\ell}
 \le C_\ell\|u\|_\ell^{\ell}.
\]
Starting with $u\in L^2(\R^N)$ and iterating while
$\ell<2N/\alpha$ yields
\[
 u\in L^\ell(\R^N)
 \qquad\text{for }2\le\ell<\frac N\alpha\,2_s^*.
\]
\end{proof}

\begin{lemma}\label{Lem2.6}
Let $r\in[(N+\alpha)/N,p]$. Suppose that $(u_n)$ is bounded in
$H^s(\R^N)$, $u_n\rightharpoonup u$ in $H^s(\R^N)$, and
$u_n\to u$ almost everywhere. If $v_n=u_n-u$, then
\begin{align}
 B_r(u_n)&=B_r(u)+B_r(v_n)+o(1),\label{eq2.1}\\
 B_r'(u_n)&=B_r'(u)+B_r'(v_n)+o(1)
 \quad\text{in }H^{-s}(\R^N),\label{eq2.2}
\end{align}
where
\[
 B_r'(z)[\varphi]
 =2r\int_{\R^N}(I_\alpha*|z|^r)|z|^{r-2}z\varphi\,dx.
\]
If $z_n\rightharpoonup0$ in $H^s(\R^N)$ and $z_n\to0$ almost
everywhere, then
\begin{equation}\label{eq2.3}
 B_r'(z_n)[\varphi]\longrightarrow0
 \qquad\text{for every fixed }\varphi\in H^s(\R^N).
\end{equation}
\end{lemma}

\begin{proof}
Identity \eqref{eq2.1} is the nonlocal Brezis--Lieb formula
\cite[Lemma~2.4]{MorozVanSchaftingen2013}. To prove the derivative
splitting, set
\[
 t=\frac{2N}{N+\alpha},\qquad
 \ell=rt,\qquad d=\frac{\ell}{r-1},\qquad k=\ell'.
\]
Then $2\le\ell\le2_s^*$, $H^s(\R^N)\hookrightarrow
L^\ell(\R^N)$, and, with $t'=2N/(N-\alpha)$,
\begin{equation}\label{eq2.4}
 \frac1k=\frac1{t'}+\frac1d.
\end{equation}
The generalized Brezis--Lieb lemma gives
\begin{align*}
 f_n&:=|u_n|^r-|v_n|^r-|u|^r\longrightarrow0
 &&\text{in }L^t(\R^N),\\
 g_n&:=|u_n|^{r-2}u_n-|v_n|^{r-2}v_n-|u|^{r-2}u
 \longrightarrow0&&\text{in }L^d(\R^N).
\end{align*}
The functional splitting in \eqref{eq2.1} follows from HLS and the fact that
$|v_n|^r\rightharpoonup0$ in $L^t$.

Set $F_n^v=|v_n|^r$, $F=|u|^r$,
$G_n^v=|v_n|^{r-2}v_n$, and $G=|u|^{r-2}u$. Expanding the derivative,
\begin{align*}
 &(I_\alpha*|u_n|^r)|u_n|^{r-2}u_n
 -(I_\alpha*F_n^v)G_n^v-(I_\alpha*F)G\\
 &\quad=(I_\alpha*f_n)|u_n|^{r-2}u_n
 +(I_\alpha*(F_n^v+F))g_n
 +(I_\alpha*F_n^v)G+(I_\alpha*F)G_n^v.
\end{align*}
The first two terms tend to zero in $L^k$ by HLS and the strong
convergence of $f_n$ and $g_n$. The last two terms tend to zero in $L^k$
by Lemma~\ref{Lem2.3} and \eqref{eq2.4}. Since $k=\ell'$ and
$H^s\hookrightarrow L^\ell$, this proves \eqref{eq2.2}.

For \eqref{eq2.3}, first take a bounded compactly supported test function
$\varphi$. Since $d>t$, one has $d'<t'$. Lemma~\ref{Lem2.3} gives
$I_\alpha*|z_n|^r\to0$ in $L^{d'}$ on its support, while
$|z_n|^{r-2}z_n$ is bounded in $L^d$. The conclusion follows from
H\"older's inequality. For a general fixed
$\varphi\in H^s\hookrightarrow L^\ell$, approximate it in $L^\ell$ by
bounded compactly supported functions. The error is uniform because
\[
 \|(I_\alpha*|z_n|^r)|z_n|^{r-2}z_n\|_k
 \le C\|z_n\|_\ell^{2r-1}.
\]
\end{proof}

\section{The autonomous problem}\label{Sec3}
In this section we consider the autonomous equation
\[
 (-\Delta)^s u=\lambda u
 +(I_\alpha*|u|^p)|u|^{p-2}u
 +(I_\alpha*|u|^q)|u|^{q-2}u
\quad\text{in }\R^N.
\]
Throughout this section we assume
\[
 \begin{gathered}
 N\ge2,\qquad 0<s<\min\left\{1,\frac N4\right\},
 \qquad N-4s<\alpha<N,\\
 \bar p:=\frac{N+2s+\alpha}{N}<q
 <p:=\frac{N+\alpha}{N-2s}.
 \end{gathered}
\]
For $u\in H^s(\R^N)$ and
$r\in[(N+\alpha)/N,p]$, set
\[
 A(u):=\|(-\Delta)^{s/2}u\|_2^2,
 \qquad
 B_r(u):=\int_{\R^N}(I_\alpha*|u|^r)|u|^r\,dx.
\]
Set
\[
 \gamma:=\frac{Nq-(N+\alpha)}{2sq},\qquad
 \beta:=q\gamma=\frac{Nq-(N+\alpha)}{2s}.
\]
The exponent assumptions imply
\begin{equation}\label{eq3.1}
 1<\beta<p,\qquad 0<\gamma<1.
\end{equation}
\begin{definition}\label{Def3.1}
For \(u\in H^s(\R^N)\), set
\[
 J_0(u):=\frac12A(u)-\frac1{2p}B_p(u)-\frac1{2q}B_q(u).
\]
For \(a>0\), let
\[
 S(a):=\{u\in H^s(\R^N):\|u\|_2^2=a\},
 \qquad (\theta\star u)(x):=e^{N\theta/2}u(e^\theta x),
\]
and define
\[
 P_0(u):=A(u)-B_p(u)-\gamma B_q(u),
 \qquad \mathcal P_a:=\{u\in S(a):P_0(u)=0\},
\]
as well as the autonomous ground level
\[
 m_a:=\inf_{u\in\mathcal P_a}J_0(u).
\]
\end{definition}

We retain the normalization of Lemma~\ref{Lem2.2}, namely
\begin{equation}\label{eq3.2}
 \mathcal S_\alpha
 =\inf_{u\in\mathcal D^{s,2}(\R^N)\setminus\{0\}}
 \frac{A(u)}{B_p(u)^{1/p}},
 \qquad B_p(u)\le \mathcal S_\alpha^{-p}A(u)^p.
\end{equation}
Moreover, Lemma~\ref{Lem2.1} gives
\begin{equation}\label{eq3.3}
 B_q(u)\le C_{\alpha,q}
 A(u)^\beta\|u\|_2^{2q(1-\gamma)}
 =C_{\alpha,q}a^{q(1-\gamma)}A(u)^\beta,
 \qquad u\in S(a).
\end{equation}

\subsection{Fiber maps and the ground level}

The logarithmic action in Definition~\ref{Def3.1} preserves \(S(a)\), and
\[
 A(\theta\star u)=e^{2s\theta}A(u),\qquad
 B_p(\theta\star u)=e^{2sp\theta}B_p(u),\qquad
 B_q(\theta\star u)=e^{2s\beta\theta}B_q(u).
\]
Consequently, the fiber map
\(\Phi_u(\theta):=J_0(\theta\star u)\) is
\[
 \Phi_u(\theta)
 =\frac12e^{2s\theta}A(u)
 -\frac1{2p}e^{2sp\theta}B_p(u)
 -\frac1{2q}e^{2s\beta\theta}B_q(u).
\]
Differentiation gives
\begin{align}
 \Phi_u'(\theta)
 &=sP_0(\theta\star u),\nonumber\\
 \Phi_u''(\theta)
 &=2s^2\big(e^{2s\theta}A(u)-p e^{2sp\theta}B_p(u)
 -\beta\gamma e^{2s\beta\theta}B_q(u)\big).
 \label{eq3.4}
\end{align}

\begin{lemma}\label{Lem3.2}
For every \(u\in S(a)\) there is a unique \(\theta(u)\in\R\) such that
\(\theta(u)\star u\in\mathcal P_a\). It is the unique critical point and
the strict global maximum of \(\Phi_u\).
Moreover,
\begin{equation}\label{eq3.5}
 \Phi_u''(\theta(u))
 =-2s^2\big((p-1)B_p(\theta(u)\star u)
 +\gamma(\beta-1)B_q(\theta(u)\star u)\big)<0.
\end{equation}
The map \(u\mapsto\theta(u)\) is \(C^1\), and
\[
 \mathcal F_a(u):=J_0(\theta(u)\star u)
 =\max_{\theta\in\R}J_0(\theta\star u)
\]
is \(C^1\) on \(S(a)\).
\end{lemma}

\begin{proof}
After division by \(e^{2s\theta}\), the sign of
\(P_0(\theta\star u)\) is the sign of
\[
 h_u(\theta)=A(u)-e^{2s(p-1)\theta}B_p(u)
 -\gamma e^{2s(\beta-1)\theta}B_q(u).
\]
By \eqref{eq3.1}, this function is strictly
decreasing from \(A(u)>0\) to \(-\infty\). This proves the unique root and
the sign change of the fiber derivative. Substituting \(P_0=0\) into
\eqref{eq3.4} gives
\eqref{eq3.5}. Finally, the scalar map
\((\theta,u)\mapsto P_0(\theta\star u)\) is \(C^1\), and its
\(\theta\)-derivative at the root is strictly negative. The implicit
function theorem gives the last assertions.
\end{proof}

By Lemma~\ref{Lem3.2}, the definition of \(m_a\) is equivalently
\[
 m_a=\inf_{u\in S(a)}\max_{\theta\in\R}J_0(\theta\star u).
\]

\begin{lemma}\label{Lem3.3}
Let $u\in H^s(\R^N)\setminus\{0\}$ be a real weak solution of
\begin{equation}\label{eq3.6}
 J_0'(u)=\lambda u\quad\text{in }H^{-s}(\R^N),
 \qquad \lambda<0.
\end{equation}
Then, for some $\eta>\max\{0,1-2s\}$ and $\vartheta>2s$,
\[
 u\in L^\infty(\R^N)\cap\operatorname{Lip}_{\mathrm{loc}}(\R^N)
 \cap C^\vartheta_{\mathrm{loc}}(\R^N),
\]
and
\[
 I_\alpha*|u|^r\in L^\infty(\R^N)\cap C^{0,\eta}(\R^N),
 \qquad r=p,q.
\]
In particular, \eqref{eq3.6} holds pointwise.
\end{lemma}

\begin{proof}
Put $\omega=-\lambda>0$ and
\[
 h_r(t)=|t|^{r-2}t,\qquad H_r=h_r(u),\qquad r=p,q.
\]
Then
\begin{equation}\label{eq3.7}
 (-\Delta)^s u+\omega u
 =\sum_{r\in\{p,q\}}(I_\alpha*(H_r u))H_r.
\end{equation}
For $r=p,q$, the exponent assumptions give
\[
 (r-1)\frac{2N}{\alpha}>2,
 \qquad
 (r-1)\frac{2N}{\alpha+2s}\le2_s^*.
\]
Since $u\in L^2(\R^N)\cap L^{2_s^*}(\R^N)$,
\[
 H_r\mathbf1_{\{|u|\le1\}}\in L^{2N/\alpha}(\R^N),
 \qquad
 H_r\mathbf1_{\{|u|>1\}}\in L^{2N/(\alpha+2s)}(\R^N).
\]
Thus Lemma~\ref{Lem2.5} applied to \eqref{eq3.7} yields
\begin{equation}\label{eq3.8}
 u\in L^t(\R^N)
 \qquad\text{for every }2\le t<t_*:=\frac N\alpha\,2_s^*.
\end{equation}

We next prove that the two Riesz potentials are bounded. For each
$r\in\{p,q\}$ choose
\begin{equation}\label{eq3.9}
 \max\left\{1,\frac2r\right\}<m_{r,-}<\frac N\alpha
 <m_{r,+}<\frac{t_*}{r}.
\end{equation}
These intervals are nonempty because $r>2\alpha/N$ and $r<2_s^*$.
By \eqref{eq3.8} and \eqref{eq3.9},
$|u|^r\in L^{m_{r,-}}\cap L^{m_{r,+}}$. Moreover,
\[
 I_\alpha\mathbf1_{B_1}\in L^{m_{r,+}'},
 \qquad
 I_\alpha\mathbf1_{\R^N\setminus B_1}\in L^{m_{r,-}'}.
\]
Young's inequality therefore gives
\begin{equation}\label{eq3.10}
 I_\alpha*|u|^r\in L^\infty(\R^N),
 \qquad r=p,q.
\end{equation}

Set
\[
 g(x,t)=(1-\omega)t+
 \sum_{r\in\{p,q\}}(I_\alpha*|u|^r)(x)|t|^{r-2}t.
\]
Then $u$ solves $(-\Delta)^s u+u=g(x,u)$ and, by
\eqref{eq3.10},
\[
 |g(x,t)|\le C\bigl(1+|t|^{2_s^*-1}\bigr).
\]
The Bessel-kernel iteration in
\cite[proof of Proposition~4, equations (3.25)--(3.33)]{Ambrosio2025Remarks}
therefore gives $u\in L^\infty(\R^N)$. Interpolation with $L^2$
then yields $u\in L^t(\R^N)$ for every finite $t\ge2$. The
right-hand side of \eqref{eq3.7} is bounded, and
\cite[Proposition~5(b)]{Ambrosio2025Remarks} gives
\begin{equation}\label{eq3.11}
 u\in C^{0,\theta_0}(\R^N)
 \qquad\text{for every }0<\theta_0<\min\{2s,1\}.
\end{equation}

Put $\sigma=\min\{1,q-1\}$. Then
\begin{equation}\label{eq3.12}
 \sigma>1-2s.
\end{equation}
This is immediate if $q\ge2$. If $q<2$, then
\[
 q-1>\frac{2s+\alpha}{N}>1-\frac{2s}{N}\ge1-2s.
\]
Since $\alpha>N-4s>1-2s$, choose $\ell$ such that
\begin{equation}\label{eq3.13}
 \ell>\max\left\{\frac N\alpha,\frac2q,
                  \frac{N}{\alpha+2s-1}\right\},
 \qquad
 \ell<\frac{N}{\alpha-1}\quad\text{if }\alpha>1.
\end{equation}
For any $\ell$ satisfying \eqref{eq3.13}, set
$\eta=\alpha-N/\ell$. Then
\begin{equation}\label{eq3.14}
 \max\{0,1-2s\}<\eta<1.
\end{equation}
Since $r\ell\ge q\ell>2$, the preceding $L^t$ estimates give
$|u|^r\in L^\ell$ for $r=p,q$. By
\cite[Theorem~2]{DuPlessis1955},
\begin{equation}\label{eq3.15}
 I_\alpha*|u|^r\in C^{0,\eta}(\R^N),
 \qquad r=p,q.
\end{equation}

Suppose that $0<\theta<1$ and $u\in C^{0,\theta}(\R^N)$. The map
$h_p$ is locally Lipschitz and $h_q$ is locally $\sigma$-H\"older.
Using \eqref{eq3.15}, the right-hand side of \eqref{eq3.7}
belongs to $C^{0,\delta}$ for every
$0<\delta<\min\{\eta,\sigma\theta\}$. Hence
\cite[Proposition~5(a)]{Ambrosio2025Remarks} gives
\begin{equation}\label{eq3.16}
 u\in C^{\theta^+}(\R^N)
 \quad\text{whenever }\theta^+<2s+\min\{\eta,\sigma\theta\},
 \qquad \theta^+\notin\mathbb N.
\end{equation}
By \eqref{eq3.12} and \eqref{eq3.14}, one may choose $c\in(0,1)$ such
that $c\eta>\max\{0,1-2s\}$ and
$c\sigma>\max\{0,1-2s\}$. Choose an exponent $\theta_0$ permitted by
\eqref{eq3.11} and apply \eqref{eq3.16} successively with
\[
 \theta_{j+1}=2s+\min\{c\eta,c\sigma\theta_j\}
\]
until the exponent exceeds one. This occurs after finitely many steps:
if the first term is selected, the next exponent is larger than one;
otherwise the affine recursion has fixed point $2s/(1-c\sigma)>1$.
Consequently,
\begin{equation}\label{eq3.17}
 u\in\operatorname{Lip}_{\mathrm{loc}}(\R^N).
\end{equation}
Using \eqref{eq3.15} once more, the right-hand side of
\eqref{eq3.7} is locally H\"older continuous. A final Schauder
estimate gives
\begin{equation}\label{eq3.18}
 u\in C^\vartheta_{\mathrm{loc}}(\R^N)
 \qquad\text{for some }\vartheta>2s.
\end{equation}
Finally,
\[
 \int_{\R^N}\frac{|u(x)|}{(1+|x|)^{N+2s}}\,dx<\infty.
\]
Together with \eqref{eq3.18}, this shows that $(-\Delta)^s u$ is
pointwise defined and that \eqref{eq3.7} holds pointwise.
\end{proof}

\begin{lemma}\label{Lem3.4}
Let \(u\in H^s(\R^N)\setminus\{0\}\) satisfy
\eqref{eq3.6}. Then
\begin{align}
 A(u)&=\lambda\|u\|_2^2+B_p(u)+B_q(u),
 \label{eq3.19}\\
 \frac{N-2s}{2}A(u)
 &=\frac N2\lambda\|u\|_2^2
   +\frac{N+\alpha}{2p}B_p(u)
   +\frac{N+\alpha}{2q}B_q(u),
  \label{eq3.20}
\end{align}
and
\begin{equation}\label{eq3.21}
 P_0(u)=A(u)-B_p(u)-\gamma B_q(u)=0.
\end{equation}
If \(u\in S(a)\), then \(J_0(u)\ge m_a\).
\end{lemma}

\begin{proof}
Testing \eqref{eq3.6} by \(u\) gives
\eqref{eq3.19}. Put \(\omega=-\lambda>0\), and
choose \(\chi\in C_c^\infty(\R^N)\) such that \(\chi=1\) on \(B_1\) and
\(\chi=0\) outside \(B_2\), and set
\[
 X_R(x):=\chi(x/R)x.
\]
For a compactly supported \(C^1\) vector field \(X\), introduce the two
deformation kernels
\begin{align*}
 K_X^s(x,y)
 &:=\frac{\operatorname{div}X(x)+\operatorname{div}X(y)}2
 -\frac{N+2s}{2}
  \frac{(X(x)-X(y))\cdot(x-y)}{|x-y|^2},\\
 K_X^\alpha(x,y)
 &:=\frac{\operatorname{div}X(x)+\operatorname{div}X(y)}2
 -\frac{N-\alpha}{2}
  \frac{(X(x)-X(y))\cdot(x-y)}{|x-y|^2}.
\end{align*}
The nonlocal integration-by-parts formulas
\cite[Propositions~6.3 and~6.5]{CingolaniGalloTanaka2024} read
\begin{align}
 \frac{C(N,s)}2\iint_{\R^{2N}}
 \frac{|u(x)-u(y)|^2}{|x-y|^{N+2s}}K_X^s(x,y)\,dx\,dy
 &=-\int_{\R^N}(-\Delta)^s u\,\nabla u\cdot X\,dx,
 \label{eq3.22}\\
 \iint_{\R^{2N}}I_\alpha(x-y)H(x)H(y)K_X^\alpha(x,y)\,dx\,dy
 &=-\int_{\R^N}(I_\alpha*H)\nabla H\cdot X\,dx.
 \label{eq3.23}
\end{align}
Lemma~\ref{Lem3.3} verifies the hypotheses of the cited
integration-by-parts formulas. In \eqref{eq3.23} we take
\(H=|u|^r\), which is
bounded and locally
Lipschitz,
\[
 \int_{\R^N}(I_\alpha*H)H\,dx=B_r(u)<\infty,
 \qquad
 (I_\alpha*H)|\nabla H|\in L^1_{\mathrm{loc}}(\R^N)
\]
by \eqref{eq3.10} and
\eqref{eq3.17}. The removed-diagonal errors in the proofs of
\eqref{eq3.22}--\eqref{eq3.23} are respectively
\(O(\rho^{2-2s})\) and \(O(\rho^\alpha)\), which tend to zero under the
standing assumptions.

For \(X=X_R\),
\[
 \sup_{R\ge1}
 \bigl(\|DX_R\|_\infty+
       \|\operatorname{div}X_R\|_\infty\bigr)<\infty.
\]
Consequently \(K_{X_R}^s\) and \(K_{X_R}^\alpha\) are uniformly bounded,
and, for almost every \(x\ne y\),
\[
 K_{X_R}^s(x,y)\longrightarrow\frac{N-2s}{2},
 \qquad
 K_{X_R}^\alpha(x,y)\longrightarrow\frac{N+\alpha}{2}.
\]
Dominated convergence in \eqref{eq3.22} is justified by
\[
 \frac{C(N,s)}2
 \frac{|u(x)-u(y)|^2}{|x-y|^{N+2s}}\in L^1(\R^{2N}),
\]
and in \eqref{eq3.23}, with \(H=|u|^r\), by
\[
 I_\alpha(x-y)|u(x)|^r|u(y)|^r\in L^1(\R^{2N}).
\]
It follows that
\begin{align}
 \lim_{R\to\infty}
 \int_{\R^N}(-\Delta)^s u\,\nabla u\cdot X_R\,dx
 &=\frac{2s-N}{2}A(u),
 \label{eq3.24}\\
 \lim_{R\to\infty}
 \int_{\R^N}(I_\alpha*|u|^r)h_r(u)\,
       \nabla u\cdot X_R\,dx
 &=-\frac{N+\alpha}{2r}B_r(u),
 \qquad r=p,q.
 \label{eq3.25}
\end{align}
Here \eqref{eq3.25} follows from
\[
 \nabla(|u|^r)=r h_r(u)\nabla u
\]
and \eqref{eq3.23}. Finally, ordinary integration by parts and
dominated convergence give
\begin{equation}\label{eq3.26}
 \lim_{R\to\infty}\int_{\R^N}u\nabla u\cdot X_R\,dx
 =-\frac12\lim_{R\to\infty}
    \int_{\R^N}u^2\operatorname{div}X_R\,dx
 =-\frac N2\|u\|_2^2.
\end{equation}

Multiply the pointwise equation by \(\nabla u\cdot X_R\), integrate, and
pass to the limit using
\eqref{eq3.24}--\eqref{eq3.26}.
Multiplication by \(-1\) gives the weak dilation identity
\[
 \frac{N-2s}{2}A(u)+\frac{N\omega}{2}\|u\|_2^2
 =\frac{N+\alpha}{2p}B_p(u)
  +\frac{N+\alpha}{2q}B_q(u).
\]
Since \(\omega=-\lambda\), this is
\eqref{eq3.20}.

Subtracting \eqref{eq3.20} from \(N/2\) times
\eqref{eq3.19} gives
\[
 sA(u)=sB_p(u)+s\gamma B_q(u),
\]
which is \eqref{eq3.21}. If \(u\in S(a)\),
then \(u\in\mathcal P_a\), so the definition of \(m_a\) gives
\(J_0(u)\ge m_a\).
\end{proof}

On \(\mathcal P_a\), the identity \(P_0(u)=0\) yields
\begin{align}
 J_0(u)
 &=\frac{p-1}{2p}B_p(u)+\frac{\beta-1}{2q}B_q(u),
 \label{eq3.27}\\
 &=\frac{\beta-1}{2\beta}A(u)
 +\left(\frac1{2\beta}-\frac1{2p}\right)B_p(u).
 \label{eq3.28}
\end{align}
Define also the nonnegative energy--Pohozaev combination
\[
 K(u):=J_0(u)-\frac1{2\beta}P_0(u)
 =\frac{\beta-1}{2\beta}A(u)
 +\left(\frac1{2\beta}-\frac1{2p}\right)B_p(u)\ge0.
\]

\begin{lemma}\label{Lem3.5}
For every \(a>0\), \(m_a>0\). Define the energy carried by one pure
upper-critical Hartree bubble by
\[
 c_H:=\frac{p-1}{2p}\mathcal S_\alpha^{p/(p-1)}
 =\frac{\alpha+2s}{2(N+\alpha)}
 \mathcal S_\alpha^{(N+\alpha)/(\alpha+2s)}.
\]
Then
\[
 m_a<c_H.
\]
Every minimizer of \(J_0|_{\mathcal P_a}\) is a constrained critical point
of \(J_0|_{S(a)}\).
\end{lemma}

\begin{proof}
If \(u\in\mathcal P_a\), then
\eqref{eq3.2},
\eqref{eq3.3}, and \(P_0(u)=0\) give
\[
 A(u)\le \mathcal S_\alpha^{-p}A(u)^p
 +\gamma C_{\alpha,q}a^{q(1-\gamma)}A(u)^\beta.
\]
Since \(p>\beta>1\), there is \(\underline A_a>0\) such that
\(A(u)\ge\underline A_a\) on \(\mathcal P_a\). Identity
\eqref{eq3.28} then yields \(m_a>0\).

If \(u\in\mathcal P_a\) attains \(m_a\), then
\(\mathcal F_a(v)\ge m_a=\mathcal F_a(u)\) for every \(v\in S(a)\), and
\(\theta(u)=0\). Hence \(d(\mathcal F_a|_{S(a)})(u)=0\). In the chain
rule all terms containing \(d\theta(u)\) have the factor
\(\Phi'_u(0)=sP_0(u)=0\). Therefore,
\(d(\mathcal F_a|_{S(a)})(u)=d(J_0|_{S(a)})(u)\), proving constrained
criticality.

It remains to prove the strict threshold. By Lemma~\ref{Lem2.2}, choose a
positive, centered, unit-scale optimizer for
\eqref{eq3.2} of the form
\[
 U(x)=C(1+|x|^2)^{-(N-2s)/2}.
\]
Its square is integrable at infinity if and only if \(N>4s\). Hence
\(U\in L^2(\R^N)\). For \(\rho>0\), let
\[
 U_\rho(x)=\rho^{(N-2s)/2}U(\rho x).
\]
Then \(A(U_\rho)=A(U)\), \(B_p(U_\rho)=B_p(U)\), and
\(\|U_\rho\|_2^2=\rho^{-2s}\|U\|_2^2\). Thus
\(\rho=(\|U\|_2^2/a)^{1/(2s)}\) places \(U_\rho\) in \(S(a)\).

For the pure critical part of the fiber, writing \(x=e^{2s\theta}\),
\[
 \max_{x>0}\left(\frac12A(U_\rho)x
 -\frac1{2p}B_p(U_\rho)x^p\right)
 =\frac{p-1}{2p}
 \left(\frac{A(U_\rho)}{B_p(U_\rho)^{1/p}}\right)^{p/(p-1)}
 =c_H.
\]
The full fiber has a finite maximizer by
Lemma~\ref{Lem3.2}, and at that maximizer the additional
term \(-e^{2s\beta\theta}B_q(U_\rho)/(2q)\) is strictly negative.
Consequently
\[
 m_a\le\max_{\theta\in\R}J_0(\theta\star U_\rho)<c_H.
\]
\end{proof}

\subsection{Mass monotonicity}

\begin{lemma}\label{Lem3.6}
If \(0<a_1<a_2\), then
\begin{equation}\label{eq3.29}
 m_{a_2}<m_{a_1}.
\end{equation}
\end{lemma}

\begin{proof}
Let \(k_0=(a_2/a_1)^{1/2}>1\), and choose
\(u_n\in\mathcal P_{a_1}\) with \(J_0(u_n)\to m_{a_1}\). From
\eqref{eq3.28}, \(A(u_n)\) is bounded above and,
as in Lemma~\ref{Lem3.5}, bounded away
from zero. We first show that
\begin{equation}\label{eq3.30}
 \liminf_{n\to\infty}B_q(u_n)>0.
\end{equation}
Otherwise, \(P_0(u_n)=0\) gives \(A(u_n)-B_p(u_n)\to0\). The sharp
inequality then implies
\(\liminf A(u_n)\ge\mathcal S_\alpha^{p/(p-1)}\), and
\eqref{eq3.27} yields
\(m_{a_1}\ge c_H\), contrary to Lemma~\ref{Lem3.5}.

For \(k\in[1,k_0]\), let \(\theta_n(k)\) be the fiber root of \(ku_n\), and
set
\[
 F_n(k):=\max_{\theta\in\R}J_0(\theta\star(ku_n)).
\]
Since \(P_0(u_n)=0\), one has \(\theta_n(1)=0\), and
\(P_0(ku_n)<0\) for \(k>1\), hence \(\theta_n(k)\le0\). These roots also
have a uniform lower bound. Their equation, divided by
\(k^2e^{2s\theta_n(k)}\), is
\begin{equation}\label{eq3.31}
 A(u_n)=k^{2p-2}e^{2s(p-1)\theta_n(k)}B_p(u_n)
 +\gamma k^{2q-2}e^{2s(\beta-1)\theta_n(k)}B_q(u_n).
\end{equation}
The left-hand side is bounded away from zero and both Hartree energies are
bounded above, so \eqref{eq3.31} precludes
\(\theta_n(k)\to-\infty\), uniformly for \(k\in[1,k_0]\).

Differentiating \(F_n\) and using the fiber equation gives
\begin{align*}
 F_n'(k)
 &=k e^{2s\theta_n(k)}A(u_n)
 -k^{2p-1}e^{2sp\theta_n(k)}B_p(u_n)
 -k^{2q-1}e^{2s\beta\theta_n(k)}B_q(u_n)\\
 &=-(1-\gamma)k^{2q-1}e^{2s\beta\theta_n(k)}B_q(u_n).
\end{align*}
By \eqref{eq3.30} and the uniform lower
bound for \(\theta_n(k)\), the right-hand side is at most a fixed negative
number for all large \(n\) and \(k\in[1,k_0]\). Hence, for some
\(d=d(a_1,a_2)>0\),
\[
 m_{a_2}\le F_n(k_0)\le F_n(1)-d
 =J_0(u_n)-d.
\]
Passing to the limit proves
\eqref{eq3.29}.
\end{proof}

\subsection{Existence and compactness of ground states}

Let \(\mathscr R(a)\) denote the nonnegative, radially symmetric,
radially nonincreasing functions in \(S(a)\). The following compactness
property is used in the autonomous minimization.

\begin{lemma}\label{Lem3.7}
If \((u_n)\subset\mathscr R(a)\) is bounded in \(H^s(\R^N)\), then, after
passing to a subsequence, \(u_n\rightharpoonup u\) in \(H^s\),
\(u_n\to u\) almost
everywhere, and
\begin{equation}\label{eq3.32}
 u_n\longrightarrow u\quad\text{in }L^r(\R^N)
 \qquad(2<r<2_s^*).
\end{equation}
In particular,
\begin{equation}\label{eq3.33}
 B_q(u_n)\to B_q(u),\qquad B_q(u_n-u)\to0.
\end{equation}
\end{lemma}

\begin{proof}
Local Rellich compactness gives convergence in \(L^r(B_R)\) for
\(2<r<2_s^*\). Radial monotonicity and the fixed \(L^2\)-mass give the
uniform radial tail estimate
\[
 |u_n(x)|^2\leq\frac{a}{|B_1|\,|x|^N}\quad(x\ne0),
\]
and hence, for every \(r>2\),
\[
 \int_{|x|>R}|u_n|^r
 \leq C_{N,r}a^{r/2}R^{-N(r-2)/2}.
\]
Combining the local convergence with this uniform tail bound proves
\eqref{eq3.32}. Since
\[
 2<\frac{2Nq}{N+\alpha}<2_s^*
\]
under the standing assumptions, HLS continuity gives
\eqref{eq3.33}.
\end{proof}

Define the autonomous ground-state set by
\[
 \mathcal G_a:=\big\{u\in\mathcal P_a:J_0(u)=m_a,
 \ (J_0|_{S(a)})'(u)=0\big\}.
\]
Define the centered positive representatives
\[
 \mathcal G_a^c:=\{w\in\mathcal G_a:w>0,
 \ w\text{ is radially nonincreasing about the origin}\}.
\]

\begin{lemma}\label{Lem3.8}
For every \(a>0\), \(\mathcal G_a^c\) is nonempty and compact in
\(H^s(\R^N)\). Moreover,
\begin{equation}\label{eq3.34}
 \mathcal G_a
 =\{\sigma T_z w:\sigma\in\{-1,1\},\ z\in\R^N,\
                    w\in\mathcal G_a^c\}.
\end{equation}
\end{lemma}

\begin{proof}
Choose \(\widetilde u_n\in\mathcal P_a\) with
\(J_0(\widetilde u_n)\to m_a\), put \(v_n=|\widetilde u_n|^*\), and let
\[
 u_n:=\theta(v_n)\star v_n\in\mathcal P_a.
\]
The fractional modulus inequality and the Riesz rearrangement inequality
give, for every \(\theta\in\R\),
\[
 J_0(\theta\star v_n)\le J_0(\theta\star\widetilde u_n).
\]
Consequently,
\[
 m_a\le J_0(u_n)
 =\max_\theta J_0(\theta\star v_n)
 \le J_0(\widetilde u_n)\longrightarrow m_a.
\]
Thus \((u_n)\subset\mathscr R(a)\) is a bounded minimizing sequence.
As in \eqref{eq3.30},
\[
 \liminf_{n\to\infty} B_q(u_n)>0.
\]
Otherwise, \(P_0(u_n)=0\), the sharp critical inequality, and
\(m_a<c_H\) give a contradiction.

By Lemma~\ref{Lem3.7}, after passing to a subsequence,
\[
 u_n\rightharpoonup u\text{ in }H^s,\qquad
 u_n\to u\quad\text{almost everywhere},\qquad B_q(u_n)\to B_q(u)>0.
\]
In particular, \(u\ne0\). Set \(r_n=u_n-u\) and \(b=\|u\|_2^2\in(0,a]\).
Hilbert splitting, Lemma~\ref{Lem2.6}, and
\eqref{eq3.33} yield
\begin{align}
 m_a&=K(u)+K(r_n)+o(1),\label{eq3.35}\\
 0&=P_0(u)+A(r_n)-B_p(r_n)+o(1).
 \label{eq3.36}
\end{align}
Pass to a further subsequence so that \(A(r_n)\to L_A\) and
\(B_p(r_n)\to L_p\).

If \(P_0(u)>0\), then \eqref{eq3.36} gives
\(L_p>L_A\). In particular \(L_A>0\), and the sharp inequality implies
\[
 L_A\ge\mathcal S_\alpha L_p^{1/p}
 >\mathcal S_\alpha L_A^{1/p},
 \qquad
 L_A>\mathcal S_\alpha^{p/(p-1)}.
\]
Since the two coefficients in \(K\) are positive and their sum is
\((p-1)/(2p)\), \eqref{eq3.35} would give
\(m_a>c_H\), a contradiction.

Hence \(P_0(u)\le0\). The fiber root satisfies \(\theta(u)\le0\), and
\(\theta\mapsto K(\theta\star u)\) is strictly increasing. Therefore
\begin{equation}\label{eq3.37}
 K(u)\ge K(\theta(u)\star u)
 =J_0(\theta(u)\star u)\ge m_b.
\end{equation}
If \(b<a\), Lemma~\ref{Lem3.6} gives \(m_b>m_a\), contradicting
\eqref{eq3.35}. Thus \(b=a\), so \(u_n\to u\) in
\(L^2\). Equations \eqref{eq3.35} and
\eqref{eq3.37} force \(K(r_n)\to0\), hence
\(A(r_n)\to0\). Thus \(u_n\to u\) strongly in \(H^s\),
\(P_0(u)=0\), and \(J_0(u)=m_a\).
By Lemma~\ref{Lem3.5}, \(u\) is a constrained critical point.
Testing its Euler equation and using
\(P_0(u)=0\) gives
\[
 \lambda_u a=-(1-\gamma)B_q(u)<0.
\]
The nontrivial nonnegative solution \(u\) is strictly positive.
Lemma~\ref{Lem3.3} shows that its equation holds
pointwise. If \(u(x_0)=0\),
then the right-hand side and the zero-order term vanish at \(x_0\), while
\[
 (-\Delta)^s u(x_0)
 =-C(N,s)\int_{\R^N}
       \frac{u(y)}{|x_0-y|^{N+2s}}\,dy<0,
\]
a contradiction. Hence \(u\in\mathcal G_a^c\).

To prove compactness, let \((u_n)\subset\mathcal G_a^c\). The coercive
identities bound \((u_n)\) in \(H^s\), and
\eqref{eq3.30} gives a uniform positive lower
bound for \(B_q(u_n)\). Lemma~\ref{Lem3.7} and the preceding splitting
argument then yield a strongly convergent subsequence. Its limit is a
centered positive ground state, so \(\mathcal G_a^c\) is compact.

We finally identify every ground state. Let \(w\in\mathcal G_a\) and
put \(v=|w|^*\). Again
\[
 A(v)\le A(w),\qquad B_r(v)\ge B_r(w)
 \quad(r\in\{p,q\}).
\]
The defining lower bound and the rearrangement inequalities give
\[
 m_a\le\max_{\theta\in\R}J_0(\theta\star v)
 \le\max_{\theta\in\R}J_0(\theta\star w)=m_a.
\]
Thus equality at the maximizing scale forces equality in every displayed
rearrangement inequality. From the
double-integral representation of the fractional seminorm,
\[
 |c-d|^2-(|c|-|d|)^2=2(|c||d|-cd),
\]
equality in the modulus inequality forces a real \(w\) to have a constant
sign. To apply the strict Riesz rearrangement theorem, set
\(f=|w|^p\). Since \(p>2\) and
\[
 \frac{2Np}{N+\alpha}=2_s^*,
\]
we have
\[
 f\in L^1(\R^N)\cap L^{2N/(N+\alpha)}(\R^N),
 \qquad B_p(w)<\infty.
\]
In particular, every positive superlevel set of the nonzero function \(f\)
has finite measure, so \(f\) vanishes at infinity in the sense required
by the equality classification.
The kernel \(I_\alpha(x)=A_{N,\alpha}|x|^{\alpha-N}\) is strictly
radially decreasing. Thus all the hypotheses of the strict equality
case \cite[Theorem~3.9]{MR1817225} are satisfied. Equality
\(B_p(v)=B_p(w)\) implies that \(f\), and hence \(|w|\), agrees almost
everywhere with a translate of its radial nonincreasing rearrangement.
After changing its sign if necessary, the preceding positivity argument
gives a strictly positive representative. This proves
\eqref{eq3.34}.
\end{proof}

\begin{lemma}\label{Lem3.9}
For every \(a>0\), there exists \(\kappa_a>0\) such that the Lagrange
multiplier of every \(w\in\mathcal G_a\) satisfies
\[
 \lambda_w\le-\kappa_a<0.
\]
\end{lemma}

\begin{proof}
Testing the Euler equation by \(w\) and using \(P_0(w)=0\) gives
\[
 \lambda_w a=-(1-\gamma)B_q(w)<0.
\]
By \eqref{eq3.34}, translation invariance of
\(B_q\), and compactness of \(\mathcal G_a^c\),
\[
b_a:=\min_{w\in\mathcal G_a^c}B_q(w)>0.
\]
Thus \(\kappa_a=(1-\gamma)b_a/a\) proves the assertion.
\end{proof}

\section{The nonautonomous problem}\label{Sec4}
Throughout this section we retain the standing assumptions
\[
 N\ge2,\qquad 0<s<\min\left\{1,\frac N4\right\},
 \qquad N-4s<\alpha<N,
 \qquad \frac{N+2s+\alpha}{N}<q<p:=\frac{N+\alpha}{N-2s},
\]
and the potential condition \eqref{eq1.2}. We use the autonomous notation
\[
 A(u)=\|(-\Delta)^{s/2}u\|_2^2,\qquad
 B_r(u)=\int_{\R^N}(I_\alpha*|u|^r)|u|^r\,dx\quad(r\in\{q,p\}),
\]
\[
 J_0(u)=\frac12A(u)-\frac1{2p}B_p(u)-\frac1{2q}B_q(u),
 \qquad Q_\varepsilon(u)=\int_{\R^N}V(\varepsilon x)|u|^2\,dx,
\]
\[
 J_\varepsilon(u):=J_0(u)+\frac12Q_\varepsilon(u).
\]
Constrained critical points of \(J_\varepsilon\) on \(S(a)\) are the real
weak solutions of \eqref{eq1.1}.

\subsection{Simultaneous cutoff}

Let \(\mathcal G_a^c\) be the compact set of centered positive autonomous
ground states obtained in Lemma~\ref{Lem3.8}, and
put
\[
 R_a:=\max_{w\in\mathcal G_a^c}\|w\|_{H^s}.
\]
Choose \(R_a<R_0<R_1\) and a nonincreasing
\(\tau\in C^\infty([0,\infty);[0,1])\) such that
\[
 \tau=1\quad\text{on }[0,R_0],\qquad
 \tau=0\quad\text{on }[R_1,\infty).
\]
On \(S(a)\), define
\begin{equation}\label{eq4.1}
 \begin{split}
 J_{\varepsilon,T}(u)
 &=\frac12A(u)+\frac12Q_\varepsilon(u)-\tau(\|u\|_{H^s})
 \left(\frac1{2p}B_p(u)+\frac1{2q}B_q(u)\right),
 \end{split}
\end{equation}
where \(T=(R_0,R_1,\tau)\). In \eqref{eq4.1}, the same cutoff
multiplies both Hartree terms. If
\(
 \mathcal N(u)=B_p(u)/(2p)+B_q(u)/(2q)
\), then
\[
 \begin{split}
 \langle J'_{\varepsilon,T}(u),\varphi\rangle
 &=\langle u,\varphi\rangle_{\dot H^s}
   +\int_{\R^N}V(\varepsilon x)u\varphi\,dx
   -\tau(\|u\|_{H^s})\langle\mathcal N'(u),\varphi\rangle\\
 &\quad-\frac{\tau'(\|u\|_{H^s})}{\|u\|_{H^s}}
       \mathcal N(u)\langle u,\varphi\rangle_{H^s}.
 \end{split}
\]
In particular, on \(\{\|u\|_{H^s}<R_0\}\),
\[
 J_{\varepsilon,T}(u)=J_0(u)+\frac12Q_\varepsilon(u),
 \qquad J'_{\varepsilon,T}(u)=J'_\varepsilon(u).
\]

\begin{lemma}\label{Lem4.1}
The zero set \(M\) is compact. For every
\(\sigma\in(0,V_\infty)\), there is \(R_\sigma>0\) such that
\begin{equation}\label{eq4.2}
 |y|\ge R_\sigma\quad\Longrightarrow\quad
 V(y)\ge V_\infty-\sigma.
\end{equation}
\end{lemma}

\begin{proof}
Equation \eqref{eq4.2} follows from the definition of
\(V_\infty\). Taking \(\sigma=V_\infty/2\) shows that the closed zero set
\(M\) is bounded and hence compact.
\end{proof}

\subsection{Local Palais--Smale analysis near the ground-state orbit}

Define the positive ground-state orbit by
\[
 \mathscr O_a:=\{T_z w:z\in\R^N,\ w\in\mathcal G_a^c\},
 \qquad T_z u=u(\,\cdot-z),
\]
and define the orbit distance by
\[
 d_a(u):=\operatorname{dist}_{H^s}(u,\mathscr O_a).
\]
This function is \(1\)-Lipschitz. Translation invariance and compactness of
\(\mathcal G_a^c\) imply
\[
 \omega_a(r):=\sup_{\substack{u\in S(a)\\d_a(u)\le r}}
 |J_0(u)-m_a|<\infty,
 \qquad \omega_a(r)\longrightarrow0\quad(r\downarrow0).
\]

Let $\kappa_a>0$ be the constant in Lemma~\ref{Lem3.9}. By compactness of
$\mathcal G_a^c$ in $L^2(\R^N)$, choose $R_a^*>0$ such that
\[
 \inf_{w\in\mathcal G_a^c}\int_{B_{R_a^*}}|w|^2\,dx\ge\frac{3a}{4}.
\]
Choose $\rho_a>0$ so small that
\begin{align}
& 4\rho_a<\min\{\sqrt{2a},R_0-R_a\},\label{eq4.3}\\
 &3\rho_a<\sqrt{\frac{3a}{4}}-\sqrt{\frac a2},\label{eq4.4}\\
 &\omega_a(2\rho_a)<\frac1{32}
 \min\{V_\infty a,\kappa_a a\},\label{eq4.5}\\
 &\sup_{\substack{u\in S(a),\ z\in\R^N,\ w\in\mathcal G_a^c\\
        \|u-T_z w\|_{H^s}\le3\rho_a}}
 \big|A(u)-B_p(u)-B_q(u)-a\lambda_w\big|
 <\frac18\kappa_a a.\label{eq4.6}
\end{align}
We also require
\begin{equation}\label{eq4.7}
 C_a\big[(4\rho_a)^{2p-2}+(4\rho_a)^{2q-2}\big]
 <\frac12c_\lambda,
 \qquad c_\lambda:=\min\{1,\kappa_a/4\},
\end{equation}
where $C_a$ is the HLS--Sobolev constant on $S(a)$. The uniform choice in
\eqref{eq4.6} follows from translation invariance, compactness
of $\mathcal G_a^c$, the testing identity, and continuity of $A$, $B_p$,
and $B_q$. Condition \eqref{eq4.7} is possible because
$p,q>1$. The $L^2$-distance between the positive and negative ground-state
orbits is at least $\sqrt{2a}$.

For \(\rho>0\), set
\[
 \mathscr U(\rho):=\{u\in S(a):d_a(u)<\rho\}.
\]
If \(u\in\overline{\mathscr U(2\rho_a)}\), the orbit approximation and
\eqref{eq4.3} give
\[
 \|u\|_{H^s}\le R_a+2\rho_a<R_0.
\]
Consequently \(\tau(\|u\|_{H^s})=1\),
\(\tau'(\|u\|_{H^s})=0\), and
\begin{equation}\label{eq4.8}
 0\le Q_\varepsilon(u)
 =2\bigl(J_{\varepsilon,T}(u)-J_0(u)\bigr)
 \le2\bigl(|J_{\varepsilon,T}(u)-m_a|+\omega_a(2\rho_a)\bigr).
\end{equation}

\begin{lemma}\label{Lem4.2}
Let \(\varepsilon_n\to0\) and \(u_n\in S(a)\) satisfy
\(d_a(u_n)\le2\rho_a\) and
\begin{equation}\label{eq4.9}
 J_{\varepsilon_n,T}(u_n)\to m_a,
 \qquad
 \|(J_{\varepsilon_n,T}|_{S(a)})'(u_n)\|\to0.
\end{equation}
Then, after passing to a subsequence, there are \(z_n\in\R^N\) and
\(w\in\mathcal G_a^c\) such that
\begin{equation}\label{eq4.10}
 T_{-z_n}u_n\to w\quad\text{in }H^s(\R^N),
 \qquad \operatorname{dist}(\varepsilon_n z_n,M)\to0.
\end{equation}
The corresponding constrained multipliers satisfy
\begin{equation}\label{eq4.11}
 \lambda_n\to\lambda_w\le-\kappa_a<0.
\end{equation}
\end{lemma}

\begin{proof}
Choose \(z_n\in\R^N\) and
\(w_n\in\mathcal G_a^c\) so that
\begin{equation}\label{eq4.12}
 \|u_n-T_{z_n}w_n\|_{H^s}\le2\rho_a+o(1).
\end{equation}
After passing to a subsequence, \(w_n\to w_0\) in \(H^s\).
By \eqref{eq4.4} and \eqref{eq4.12},
\[
 \begin{aligned}
 \|u_n\|_{L^2(B_{R_a^*}(z_n))}
 &\ge \|w_n\|_{L^2(B_{R_a^*})}
      -\|u_n-T_{z_n}w_n\|_2\\
 &\ge \sqrt{\frac{3a}{4}}-2\rho_a-o(1)
 >\sqrt{\frac a2}
 \end{aligned}
\]
for all large $n$. The sequence $y_n=\varepsilon_n z_n$ is bounded.
Otherwise, a subsequence would satisfy $|y_n|\to\infty$.
Lemma~\ref{Lem4.1} and \eqref{eq4.2} would then give
\begin{equation}\label{eq4.13}
 \liminf_{n\to\infty} Q_{\varepsilon_n}(u_n)\ge\frac14V_\infty a,
\end{equation}
contrary to \eqref{eq4.8},
\eqref{eq4.5}, and \eqref{eq4.9}.
After passing to a subsequence, \(y_n\to y\).

Put \(v_n=T_{-z_n}u_n\). Then \(v_n\rightharpoonup v\) in \(H^s\) and
\(v_n\to v\) almost everywhere. Moreover,
\(\|v-w_0\|_{H^s}\le2\rho_a\) by weak lower semicontinuity and
\eqref{eq4.12}. Thus \(v\ne0\). Otherwise,
\(\sqrt a\le\|w_0\|_{H^s}\le2\rho_a\), contrary to
\eqref{eq4.3}. Lemma~\ref{Lem2.4} gives, for every fixed
$\varphi\in H^s(\R^N)$,
\[
 \|[V(y_n+\varepsilon_n\,\cdot)-V(y)]\varphi\|_2\to0.
\]
Choose multipliers \(\lambda_n\) so that
\begin{equation}\label{eq4.14}
 J'_{\varepsilon_n,T}(u_n)-\lambda_n u_n\to0
 \quad\text{in }H^{-s}.
\end{equation}
Testing by \(u_n\), and using \eqref{eq4.8},
\eqref{eq4.6}, \eqref{eq4.5}, and
\eqref{eq4.9}, yields
\[
 \lambda_n a=A(u_n)+Q_{\varepsilon_n}(u_n)-B_p(u_n)-B_q(u_n)+o(1),
 \qquad \lambda_n\le-\frac12\kappa_a+o(1).
\]
The same identity bounds \((\lambda_n)\). After passing to a subsequence,
let \(\lambda_n\to\lambda<0\). The derivative splitting in
Lemma~\ref{Lem2.6}, specifically
\eqref{eq2.3} for the weakly null remainder,
and Lemma~\ref{Lem2.4} permit passage to the limit against every fixed
test function in the translated equation. Thus
\begin{equation}\label{eq4.15}
 J_0'(v)=(\lambda-V(y))v.
\end{equation}

Set \(r_n=v_n-v\). The Hilbert splitting around \(w_n\to w_0\) gives
\(\limsup\|r_n\|_{H^s}\le4\rho_a\). Subtracting
\eqref{eq4.15} from \eqref{eq4.14}, testing by
\(r_n\), and using Lemma~\ref{Lem2.6} gives
\[
 A(r_n)+\int_{\R^N}V(y_n+\varepsilon_n x)|r_n|^2\,dx
 -\lambda_n\|r_n\|_2^2
 =B_p(r_n)+B_q(r_n)+o(1).
\]
The cross terms are $o(1)$ by Lemma~\ref{Lem2.4}, weak
\(L^2\)-convergence, and boundedness of \(V\). Since
\[
 B_p(r_n)+B_q(r_n)
 \le C_a\bigl(\|r_n\|_{H^s}^{2p-2}
              +\|r_n\|_{H^s}^{2q-2}\bigr)\|r_n\|_{H^s}^2,
\]
for all large $n$ the left-hand side is at least
$c_\lambda\|r_n\|_{H^s}^2$. Thus \eqref{eq4.7} implies
$r_n\to0$ in $H^s$.

Strong convergence gives \(\|v\|_2^2=a\). Since
$v_n\to v$ in $L^2(\R^N)$, Lemma~\ref{Lem2.4} gives
\[
 Q_{\varepsilon_n}(u_n)
 =\int_{\R^N}V(y_n+\varepsilon_n x)|v_n(x)|^2\,dx
 \longrightarrow V(y)a.
\]
Together with the nonlinear and kinetic strong convergence, this gives
\begin{equation}\label{eq4.16}
 m_a=J_0(v)+\frac a2V(y).
\end{equation}
The multiplier in \eqref{eq4.15} is negative, so
Lemma~\ref{Lem3.4} gives
\(P_0(v)=0\), hence \(J_0(v)\ge m_a\). Equation
\eqref{eq4.16} and \(V\ge0\) force
\(V(y)=0\) and \(J_0(v)=m_a\). Therefore \(v\) is a ground state. The
orbit representation \eqref{eq3.34} and
\(4\rho_a<\sqrt{2a}\) exclude the negative orbit. Absorbing the remaining
bounded translation into \(z_n\) proves
\eqref{eq4.10}. Since $V(y)=0$, comparing \eqref{eq4.15} with
the Euler equation of $v$ gives $\lambda=\lambda_v$. This proves
\eqref{eq4.11}.
\end{proof}

\begin{lemma}\label{Lem4.3}
There exist \(\eta_a,\varepsilon_a,\kappa_{\rm gap}>0\) such that
\begin{equation}\label{eq4.17}
 2\bigl(\eta_a+\omega_a(2\rho_a)\bigr)
 <\frac14\min\{V_\infty a,\kappa_a a\},
\end{equation}
and the following assertions hold.
\begin{enumerate}
\item[(i)] For every fixed \(0<\varepsilon<\varepsilon_a\),
\(J_{\varepsilon,T}|_{S(a)}\) satisfies the Palais--Smale condition for
sequences in
\begin{equation}\label{eq4.18}
 \overline{\mathscr U(2\rho_a)}
 \cap J_{\varepsilon,T}^{-1}([m_a-\eta_a,m_a+\eta_a]).
\end{equation}
Their multipliers satisfy
\begin{equation}\label{eq4.19}
 \lambda_n\le-\frac12\kappa_a+o(1).
\end{equation}

\item[(ii)] For every \(0<\varepsilon<\varepsilon_a\),
\[
 \|(J_{\varepsilon,T}|_{S(a)})'(u)\|\ge\kappa_{\rm gap}
\]
whenever
\[
 \rho_a\le d_a(u)\le2\rho_a,
 \qquad |J_{\varepsilon,T}(u)-m_a|\le\eta_a.
\]
\end{enumerate}
\end{lemma}

\begin{proof}
Choose \(\eta_*>0\) so small that
\[
 2\bigl(\eta_*+\omega_a(2\rho_a)\bigr)
 <\frac14\min\{V_\infty a,\kappa_a a\}.
\]
We first prove the fixed-parameter Palais--Smale property on the same
domain as \eqref{eq4.18}, with \(\eta_*\) in place of
\(\eta_a\). Let \((u_n)\) be such a fixed-\(\varepsilon\) PS sequence and
choose \(z_n,w_n\) as in
\eqref{eq4.12}. If \(|z_n|\to\infty\), then
\eqref{eq4.2} and the local-mass estimate used in
\eqref{eq4.13} give
\(
 \liminf Q_\varepsilon(u_n)\ge V_\infty a/4
\), contradicting \eqref{eq4.8} and the defining inequality
for \(\eta_*\). Thus \((z_n)\) is bounded. After passing to a subsequence,
\(z_n\to z\), \(w_n\to w_0\), and
\(u_n\rightharpoonup u\ne0\) in \(H^s(\R^N)\). By the local Rellich
compactness theorem and a diagonal argument,
\[
 u_n\longrightarrow u\quad\text{in }L^r_{\mathrm{loc}}(\R^N)
 \ \text{for every }2\le r<2_s^*,
 \qquad u_n\longrightarrow u\quad\text{almost everywhere in }\R^N.
\]

Choose the constrained multipliers so that, using the inactivity of the
cutoff throughout the tube,
\begin{equation}\label{eq4.20}
 J_\varepsilon'(u_n)-\lambda_n u_n\longrightarrow0
 \qquad\text{in }H^{-s}(\R^N).
\end{equation}
Testing by \(u_n\), using \eqref{eq4.8} and the
three-\(\rho_a\) margin in \eqref{eq4.6}, gives
\[
 \lambda_n a=A(u_n)+Q_\varepsilon(u_n)-B_p(u_n)-B_q(u_n)+o(1),
 \qquad
 \lambda_n\le-\frac12\kappa_a+o(1).
\]
Thus \((\lambda_n)\) is bounded. After passing to a subsequence, let
\(\lambda_n\to\lambda<0\). Lemma~\ref{Lem2.6}, in particular the
weak-null conclusion \eqref{eq2.3}, permits
passage to the limit in \eqref{eq4.20}. Hence
\begin{equation}\label{eq4.21}
 J_\varepsilon'(u)=\lambda u
 \qquad\text{in }H^{-s}(\R^N).
\end{equation}

Set \(r_n=u_n-u\). The bounded centers and the orbit approximation give
\(\limsup_{n\to\infty}\|r_n\|_{H^s}\le4\rho_a\). Subtracting
\eqref{eq4.21} from \eqref{eq4.20}, testing
by \(r_n\), and using the full splitting in Lemma~\ref{Lem2.6} gives
\[
 A(r_n)+\int_{\R^N}V(\varepsilon x)|r_n|^2\,dx
 -\lambda_n\|r_n\|_2^2
 =B_p(r_n)+B_q(r_n)+o(1).
\]
The potential cross term is \(o(1)\) because
\(V(\varepsilon\,\cdot)u\in L^2\) and \(r_n\rightharpoonup0\) in \(L^2\).
For all large $n$, the nonnegative potential and
$-\lambda_n\ge\kappa_a/4$ bound the left-hand side from below by
$c_\lambda\|r_n\|_{H^s}^2$. Condition \eqref{eq4.7} then
gives $r_n\to0$ in $H^s$. This proves the fixed-parameter local
Palais--Smale property and \eqref{eq4.19}.

We prove the uniform annular bound. Suppose, for contradiction,
that no triple
\((\eta_{\rm gap},\varepsilon_{\rm gap},\kappa_{\rm gap})\) yields the
asserted estimate. Then for every \(n\) there are
\(0<\varepsilon_n<1/n\) and \(u_n\in S(a)\) such that
\[
 \rho_a\le d_a(u_n)\le2\rho_a,\qquad
 |J_{\varepsilon_n,T}(u_n)-m_a|
 \le\min\{\eta_*,1/n\},\qquad
 \|(J_{\varepsilon_n,T}|_{S(a)})'(u_n)\|\le\frac1n.
\]
Lemma~\ref{Lem4.2} would give \(d_a(u_n)\to0\), contradicting
\(d_a(u_n)\ge\rho_a\). Hence such a triple exists. Set
\[
 \eta_a:=\min\{\eta_*,\eta_{\rm gap}\},
 \qquad \varepsilon_a:=\varepsilon_{\rm gap}.
\]
Then \eqref{eq4.17}, the fixed-parameter result, and the
annular estimate all hold with these constants.
\end{proof}

\begin{lemma}\label{Lem4.4}
For every \(u\in\overline{\mathscr U(2\rho_a)}\),
\[
 \|u\|_{H^s}<R_0,\qquad \tau(\|u\|_{H^s})=1,
 \qquad \tau'(\|u\|_{H^s})=0.
\]
Consequently every constrained critical point in this tube solves the
original, untruncated equation. If in addition
\(0<\varepsilon<\varepsilon_a\) and
\(|J_{\varepsilon,T}(u)-m_a|\le\eta_a\), then its multiplier satisfies
\[
 \lambda_u\le-\frac12\kappa_a<0.
\]
\end{lemma}

\begin{proof}
For \(u\in\overline{\mathscr U(2\rho_a)}\),
\eqref{eq4.3} gives
\(\|u\|_{H^s}\le R_a+2\rho_a<R_0\), which proves the cutoff assertions.
The multiplier estimate follows by applying
Lemma~\ref{Lem4.3}(i) to the constant Palais--Smale sequence \(u_n=u\).
\end{proof}

\section{Existence and concentration}\label{Sec5}

Fix \(y_0\in M\) and \(w\in\mathcal G_a^c\). For \(t>0\), set
\[
 w_t(x):=t^{N/2}w(tx)\in S(a).
\]
By Lemma~\ref{Lem3.2},
\begin{equation}\label{eq5.1}
 J_0(w_t)<m_a\quad(t\ne1),\qquad
 P_0(w_t)>0\quad(0<t<1),\qquad
 P_0(w_t)<0\quad(t>1).
\end{equation}

Choose \(t_-<1<t_+\) sufficiently close to \(1\) that
\begin{equation}\label{eq5.2}
 \max_{t\in[t_-,t_+]}\|w_t-w\|_{H^s}<\frac14\rho_a.
\end{equation}
After this choice, fix \(\eta_0>0\) such that
\begin{equation}\label{eq5.3}
 4\eta_0<\min\{m_a-J_0(w_{t_-}),
                 m_a-J_0(w_{t_+}),\eta_a\}.
\end{equation}
\begin{lemma}\label{Lem5.1}
For every \(0<\varepsilon<\varepsilon_a\), define
\[
 \Gamma_\varepsilon(t):=T_{y_0/\varepsilon}w_t,
 \qquad t\in[t_-,t_+],
\]
and set
\[
 h_\varepsilon:=\varepsilon+\frac12
 \sup_{t\in[t_-,t_+]}Q_\varepsilon(\Gamma_\varepsilon(t)).
\]
Then
\[
 \Gamma_\varepsilon([t_-,t_+])
 \subset\mathscr U(\rho_a/4),\qquad h_\varepsilon\longrightarrow0.
\]
After decreasing \(\varepsilon_a\) if necessary, the following estimates
hold for every \(0<\varepsilon<\varepsilon_a\):
\begin{align}
 \max_{t\in[t_-,t_+]}
 J_{\varepsilon,T}(\Gamma_\varepsilon(t))
 &\le m_a+h_\varepsilon,
 \label{eq5.4}\\
 J_{\varepsilon,T}(\Gamma_\varepsilon(t_-)),\quad
 J_{\varepsilon,T}(\Gamma_\varepsilon(t_+))
 &\le m_a-3\eta_0,
 \label{eq5.5}
\end{align}
and
\begin{equation}\label{eq5.6}
 h_\varepsilon<\min\{\eta_0,\eta_a,
                         \tfrac12\kappa_{\rm gap}\rho_a\}.
\end{equation}
\end{lemma}

\begin{proof}
The set $\{w_t:t\in[t_-,t_+]\}$ is compact in
$L^2(\R^N)$. Since $V(y_0)=0$, Lemma~\ref{Lem2.4} gives
\begin{equation}\label{eq5.7}
 \sup_{t\in[t_-,t_+]}
 Q_\varepsilon(\Gamma_\varepsilon(t))\longrightarrow0.
\end{equation}
The inclusion in the tube follows from \eqref{eq5.2}, so the
cutoff is inactive on the whole path. Hence
\[
 J_{\varepsilon,T}(\Gamma_\varepsilon(t))
 =J_0(w_t)+\frac12Q_\varepsilon(\Gamma_\varepsilon(t)).
\]
Now \eqref{eq5.1}, \eqref{eq5.3}, and
\eqref{eq5.7} yield
\eqref{eq5.4}--\eqref{eq5.5}.
Since \(h_\varepsilon\to0\), the constant \(\varepsilon_a\) may be
decreased so that all these conclusions and
\eqref{eq5.6} hold throughout
\((0,\varepsilon_a)\). This reduction preserves
Lemmas~\ref{Lem4.3} and~\ref{Lem4.4}.
\end{proof}

We use a local deformation because the constrained energy is unbounded
from below on \(S(a)\).

\begin{definition}
Let \(e_\varepsilon^\pm=\Gamma_\varepsilon(t_\pm)\) and
\[
 \mathfrak G_\varepsilon
 :=\{\gamma\in C([0,1],S(a)):
       \gamma(0)=e_\varepsilon^-,\ \gamma(1)=e_\varepsilon^+\}.
\]
Define
\[
 c_\varepsilon
 :=\inf_{\gamma\in\mathfrak G_\varepsilon}
      \max_{t\in[0,1]}J_{\varepsilon,T}(\gamma(t)).
\]
\end{definition}

\begin{lemma}\label{Lem5.3}
For every \(0<\varepsilon<\varepsilon_a\),
\begin{equation}\label{eq5.8}
 m_a\le c_\varepsilon\le m_a+h_\varepsilon,
\end{equation}
and \(J_{\varepsilon,T}|_{S(a)}\) has a critical point \(u_\varepsilon\)
such that
\begin{equation}\label{eq5.9}
 d_a(u_\varepsilon)<\rho_a,
 \qquad
 m_a\le J_{\varepsilon,T}(u_\varepsilon)
 \le m_a+h_\varepsilon.
\end{equation}
\end{lemma}

\begin{proof}
Every \(\gamma\in\mathfrak G_\varepsilon\) crosses \(\{P_0=0\}\) by
\eqref{eq5.1}. Since
\[
 J_{\varepsilon,T}(u)
 =J_0(u)+\frac12Q_\varepsilon(u)
 +(1-\tau(\|u\|_{H^s}))
 \left(\frac1{2p}B_p(u)+\frac1{2q}B_q(u)\right)
 \ge J_0(u),
\]
and \(J_0(u)\ge m_a\) on \(\mathcal P_a=\{P_0=0\}\), the lower bound in
\eqref{eq5.8} follows. For the upper bound, the reparametrized
trial path
\[
 s\longmapsto
 \Gamma_\varepsilon\bigl(t_-+s(t_+-t_-)\bigr),
 \qquad 0\le s\le1,
\]
belongs to \(\mathfrak G_\varepsilon\), and
Lemma~\ref{Lem5.1} applies.

Assume, to the contrary, that there is no critical point in
\begin{equation}\label{eq5.10}
 \overline{\mathscr U(2\rho_a)}\cap
 J_{\varepsilon,T}^{-1}([m_a,m_a+h_\varepsilon]).
\end{equation}
By the fixed-parameter local PS property in
Lemma~\ref{Lem4.3}(i), the absence of critical points implies that there
is \(\sigma_\varepsilon>0\) such that
\[
 \|(J_{\varepsilon,T}|_{S(a)})'(u)\|\ge\sigma_\varepsilon
\]
on the closed set \eqref{eq5.10}.

Let $K_\varepsilon$ be the closed set in
\eqref{eq5.10}. The sphere $S(a)$ is a complete smooth
Hilbert submanifold of $H^s(\R^N)$. A locally finite
pseudo-gradient construction gives a bounded, locally Lipschitz tangent
field $Y_0$ on an open neighborhood $\mathcal O_\varepsilon$ of
$K_\varepsilon$ such that
\[
 \|Y_0(u)\|_{H^s}\le1,
 \qquad J'_{\varepsilon,T}(u)[Y_0(u)]\ge0,
\]
and, after shrinking \(\mathcal O_\varepsilon\) if necessary,
\[
 J'_{\varepsilon,T}(u)[Y_0(u)]
 \ge\frac12\|(J_{\varepsilon,T}|_{S(a)})'(u)\|.
\]
Choose a locally Lipschitz tube--energy cutoff
\(\zeta_\varepsilon:S(a)\to[0,1]\), supported in
\(\mathcal O_\varepsilon\), and equal to one on \(K_\varepsilon\).
In particular, it equals one on the entire high-energy annulus
\[
 \rho_a\le d_a(u)\le2\rho_a,\qquad
 m_a\le J_{\varepsilon,T}(u)\le m_a+h_\varepsilon.
\]
The field $\zeta_\varepsilon Y_0$, extended by zero outside
$\mathcal O_\varepsilon$, is globally bounded and locally Lipschitz.
Multiply it by a smooth energy cutoff that equals one on $[m_a,\infty)$
and zero on $(-\infty,m_a-2\eta_0]$, and denote the resulting tangent
field by $Y$. Since $S(a)$ is complete and $Y$ is bounded and locally
Lipschitz, its negative flow $\Psi(s,u)$ exists for every finite time,
depends continuously on $(s,u)$, and does not increase the energy. While
the trajectory remains in \(\overline{\mathscr U(2\rho_a)}\) with
\(m_a\le J_{\varepsilon,T}(\Psi(s,u))\le m_a+h_\varepsilon\), it satisfies
\begin{equation}\label{eq5.11}
 \frac d{ds}J_{\varepsilon,T}(\Psi(s,u))
 \le-\frac12\|(J_{\varepsilon,T}|_{S(a)})'(\Psi(s,u))\|.
\end{equation}

While its energy is at least \(m_a\), a trajectory starting from the trial
arc cannot leave \(\mathscr U(2\rho_a)\). The initial arc lies in
\(\mathscr U(\rho_a/4)\). Since \(d_a\) is \(1\)-Lipschitz and the flow
speed is at most one, a trajectory exiting the outer tube must spend at
least time \(\rho_a\) while crossing
\(\rho_a\le d_a\le2\rho_a\). On this part of the trajectory
Lemma~\ref{Lem4.3}(ii) and
\eqref{eq5.11} force an energy loss of at least
\(\kappa_{\rm gap}\rho_a/2>h_\varepsilon\), contradicting the initial
upper bound \eqref{eq5.4} as long as the energy remains at least
\(m_a\). Thus no trajectory leaves the outer tube before its energy falls
below \(m_a\).

As long as a trajectory stays at or above \(m_a\), it therefore remains in
the outer tube and its energy decreases at rate at least
\(\sigma_\varepsilon/2\). For any time $T>2h_\varepsilon/\sigma_\varepsilon$, every point of
the trial arc is therefore mapped below $m_a$. The endpoints are fixed by
\eqref{eq5.5}, \eqref{eq5.6}, and the
energy cutoff. We therefore obtain a path
\(\widehat\Gamma_\varepsilon\in\mathfrak G_\varepsilon\) satisfying
\[
 \max_{t\in[0,1]} J_{\varepsilon,T}(\widehat\Gamma_\varepsilon(t))<m_a,
\]
contradicting the lower bound in \eqref{eq5.8}. Thus a critical
point exists in \eqref{eq5.10}. It cannot lie in the
annulus \(\rho_a\le d_a\le2\rho_a\) because of
Lemma~\ref{Lem4.3}(ii), so it satisfies
\eqref{eq5.9}.
\end{proof}

\begin{proof}[Proof of Theorem~\ref{Thm1.1}]
By Lemma~\ref{Lem5.3}, there exists \(u_\varepsilon\). The
cutoff is inactive by Lemma~\ref{Lem4.4}, so
\(u_\varepsilon\) is a constrained critical point of the original
functional \(J_\varepsilon\). Its multiplier is negative by
Lemma~\ref{Lem4.4}. Since the functional and the mass constraint are
even, \(-u_\varepsilon\) is a solution with the same multiplier. The two
solutions are distinct because \(a>0\).

Let \(\varepsilon_n\to0\). By \eqref{eq5.9} and
\(h_{\varepsilon_n}\to0\), the sequence \(u_{\varepsilon_n}\) satisfies
the hypotheses of Lemma~\ref{Lem4.2}. Hence, after passing to a
subsequence, there are \(z_n\in\R^N\) and \(w\in\mathcal G_a^c\) such that
\[
 u_{\varepsilon_n}(\,\cdot+z_n)\to w
 \quad\text{strongly in }H^s(\R^N),
 \qquad
 \operatorname{dist}(\varepsilon_n z_n,M)\to0,
\]
and \(\lambda_{\varepsilon_n}\to\lambda_w<0\). This is the asserted
concentration statement.
\end{proof}
\section*{Acknowledgments}

This work was supported by the National Natural Science
Foundation of China (12301145, 12161007, 12261107), the Guangxi Natural
Science Foundation (2023GXNSFAA026190), and the Yunnan Fundamental
Research Projects (202301AU070144, 202401AU070123). Y.~Aikyn was supported
by the Bolashak Government Scholarship of the Republic of Kazakhstan.
M.~Ruzhansky was supported by the FWO Odysseus 1 Grant G.0H94.18N and the
Methusalem Programme of the Ghent University Special Research Fund
(01M01021).

\medskip
\noindent\emph{Data availability:}
No data were generated or analyzed in this study.

\medskip
\noindent\emph{Conflict of interest:}
The authors declare no conflict of interest.

\medskip
\noindent\emph{AI assistance statement:}
The authors used AI models to assist with writings; the main ideas, mathematical validation, and all final checks remain the sole responsibility of the human authors.

\enlargethispage{2\baselineskip}

\end{document}